\documentclass[journal]{IEEEtran}
 \usepackage{times, mathptm, graphicx, amsmath, amsfonts, bm}

 \def\BEA{\begin{eqnarray}}
 \def\EEA{\end{eqnarray}}
 \def\BAL{\begin{align}}
 
 \def\NN{\nonumber}

 \def\EG{{\it e.g.}}
 \def\EEG{{\it E.g.}}
 \def\IE{{\it i.e.}}
 
 \def\s{,\,}
 \def\d{{\rm d}}
 \def\e{{\rm e}}
 \def\i{{\rm i}}

 \def\VEC#1{{\pmb{#1}}}

 \def\PHI{\varphi}
 \def\EPS{\varepsilon}
 \def\DOT#1{\smash{\overset{\mbox{\large .}}{#1}}}
 \def\DDOT#1{\smash{\overset{\mbox{\large ..}}{#1}}}
 \def\SEP{\,/\,}

\begin{document}

\title{Most probable failure scenario in a model\\
power grid with random power demand}

\author{Misha Stepanov and Aditya Sundarrajan \thanks{This work was 
supported by DTRA/DOD grant BRCALL06-Per3-D-2-0022 ``Network 
adaptability from WMD disruption and cascading failures''.}
 \thanks{M.~Stepanov is with Department of Mathematics, University of 
Arizona, Tucson, AZ 85721, USA (e-mail: stepanov@math.arizona.edu).}
 \thanks{A.~Sundarrajan is with Department of Electrical and Computer 
Engineering, University of Arizona, Tucson, AZ 85721, USA (e-mail: 
adityas@email.arizona.edu).}}

\maketitle

 \begin{abstract}
   We consider a simple system with a local synchronous generator and a 
load whose power consumption is a random process. The most probable 
scenario of system failure (synchronization loss) is considered, and it 
is argued that its knowledge is virtually enough to estimate the 
probability of failure per unit time. We discuss two numerical methods 
to obtain the ``optimal'' evolution leading to failure.
 \end{abstract}


\section{Introduction} \label{sec:introduction}

Whenever we use or design a system we are curious about its reliability. 
Failures may be very rare, still they could be of a tremendous 
importance, especially when the consequences are catastrophic and very 
undesirable. There is a huge amount of literature devoted to the study 
of power systems failures, for various situations and approaches see, 
\EG, \cite{1989_SCL_DC, 2007_C_DCLN, 2009_IEEETPS_SWB}. Failures could 
happen due to a hardware fault, because the realization of the system's 
random component (noise) was an especially unlucky one (the latter is a 
focus of our work), or a combination of both.

While studying such unlucky realizations, one can ask the following 
questions: (1) What is the spatio-temporal shape{\SEP}profile of the 
noise configuration that leads to a rare event{\SEP}failure? (2) Is this 
shape universal? --- If two unrelated failure events happen, will their 
pre-histories be alike?

Theoretical exploration of this subject started from the classical 
papers \cite{1964_AP_L, 1966_PR_HL, 1966_PR_ZL}. The theoretical 
approach is known under the names of optimal fluctuation or instanton 
method; it could also be viewed as a saddle-point method in functional 
space. It was originally introduced by I.~M.~Lifshitz \cite{1964_AP_L} 
for analysis of electronic spectra of disordered systems, and was later 
developed to describe rare events in various fields, \EG, quantum field 
theory \cite{1975_PLB_BPST, 1977_JETP_L}, statistical hydrodynamics 
\cite{1996_PRE_FKLM}, and power grids \cite{2011_IEEETSG_CPS, 
2011_CDC_CSPB}.

The main conjecture of the optimal fluctuation method is a positive 
answer to the second question posed above. The rationale behind this 
theoretical expectation is as follows: In order to cause such a rare 
event, the deviation of the system's random characteristics from their 
average values (\IE, the noise) should be much larger than typical. On 
the other hand statistical weight, or probability, of the noise sharply 
decays with the noise amplitude increase. On the whole variety of noise 
configurations leading to the event of interest the probability is 
highly peaked at the most probable or ``optimal'' configuration, on 
which the ``amplitude'' is not so very large, but the rare event still 
does happen due to carefully optimized shape of the noise.

Optimal fluctuation method consists in maximizing the probability of the 
random noise realization conditioned by the statement that the rare 
event of interest happens. In many cases one can come up with the 
[obviously problem specific] equations that the position of the maximum 
(\IE, optimal fluctuation) should satisfy (with the meaning of the 
equations being ``the derivative of the probability weight at the 
extremum is equal to zero''). These equations are often too complex for 
analytical solution and their numerical exploration remains the only 
feasible method of analysis.

The knowledge of optimal noise configuration can be indispensable for 
new system design, as it delivers important information on how dangerous 
and undesirable events are going to look like. Also, the probability 
[per unit time] of the rare event can be estimated.

In this work we consider a simple power system (Sec.~\ref{sec:system}) 
and derive the equations the most probable failure scenario should 
satisfy (Sec.~\ref{sec:equations}, see also Sec.~\ref{sec:HJB}). We 
present two numerical methods to solve these equations 
(Secs.~\ref{sec:IVP} and \ref{sec:BVP}), and then discuss the 
probability of failure (Sec.~\ref{sec:probability}).

\section{Equations of motion} \label{sec:system}

\begin{figure}[b]
  \begin{picture}(252,70)(0,0)
    \put(14,0){\includegraphics[width=216pt]{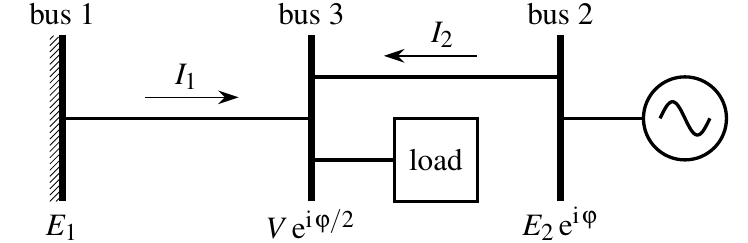}}
  \end{picture}
  \caption{Model power system. The currents $I_1$ and $I_2$, and the 
voltage $V$ are complex numbers. \label{scheme}} \end{figure}

Let us consider the following simple power system: It has 3 buses, one 
of them is a slack (or infinite) bus, another is a local synchronous 
generator bus, and the first two are connected by transmission lines to 
the 3rd one (a load bus), see Fig.~\ref{scheme}. The dynamics of the 
system is governed by
 \begin{align}
   & M \DDOT\PHI + D \DOT\PHI = P_{\rm mech} - {\rm Re} \big( E_2 \, 
\e^{\i \, \PHI} \, I_2^* \big) , \label{swing} \\
   & \i \, X_1 I_1 = E_1 - V \, \e^{\i \, \PHI / 2}, \quad \i \, X_2 
I_2 = E_2 \, \e^{\i \, \PHI} - V \e^{\i \, \PHI / 2} , \NN \\
   & V \, \e^{\i \, \PHI / 2} (I_1 + I_2)^* = S(t) = P(t) + \i \, Q(t) 
, \NN
 \end{align} where $X_1$ and $X_2$ are reactances of the lines, $P_{\rm 
mech}$ is an input of mechanical power which we assume constant, and 
$P${\SEP}$Q$ is the real{\SEP}reactive power consumed by the load. The 
eq.~\eqref{swing} is the swing equation with $M$ being the moment of 
inertia of the generator divided by a square of the number of ``poles'' 
(\IE, how many oscillations of voltage one rotation of the generator 
produces) and multiplied by angular frequency of rotation. The term $D 
\DOT\PHI$ describes damping in the generator. The dot over a variable 
denotes a time derivative, two dots mean second time derivative. The 
asterisk $\cdot^{*}$ represents the complex conjugate. By excluding 
$I_1$ and $I_2$ we get
 \begin{align}
   & M \DDOT\PHI + D \DOT\PHI = P_{\rm mech} - (E_2 / X_2) \, {\rm Im} 
\big( V^* \, \e^{\i \, \PHI / 2} \big) , \NN \\
   & S = \i \, V \Big( \frac{E_1 \, \e^{\i \, \PHI / 2}}{X_1} + 
\frac{E_2 \, \e^{-\i \, \PHI / 2}}{X_2} \Big) - \i \, |V|^2 \Big( 
\frac{1}{X_1} + \frac{1}{X_2} \Big) . \NN
 \end{align} Let us consider $X_1 = X_2 = 1$, $E_1 = E_2 = 1$, $M = 1$. 
We will also write $V = y - \i \, x$. We have
 \begin{align}
   & \DDOT\PHI + D \DOT\PHI = P_{\rm mech} - {\sf C} x - {\sf S} y , \NN 
\\
   & P = 2 {\sf C} x , \quad Q = 2 {\sf C} y - 2 (x^2 + y^2) . \NN
 \end{align} where ${\sf C} = \cos (\PHI / 2)$, ${\sf S} = \sin (\PHI / 
2)$. The boundary for the solution existence for $x$ and $y$ is given by 
the inequality $Q \le ({\sf C}^2 - P^2 / {\sf C}^2) / 2$. We will assume 
$Q = k P$ with some constant $k$ (constant power load model 
\cite{1993_IEEETPS_TaskForce}, although we will allow the power demand 
$P$ to change in time). We will use the value $k = 3 / 4$, and the 
inequality becomes $P \le {\sf C}^2 / 2$.

Let us consider as a Normal Operating Point (NOP) the steady-state (\IE, 
$\DOT\PHI = \DDOT\PHI = 0$) solution with
 \begin{align}
                  \PHI_* & = 0 , \quad {\sf C}_* = 1 , \quad {\sf S}_* = 
0 , \NN \\
   S_* = P_* + \i \, Q_* & = (4 + 3 \i) / 16 , \quad P_{\rm mech} = 1 / 
8 , \NN \\
   V_* = y_* - \i \, x_* & = (7 - \i) / 8 . \NN
 \end{align}

The dynamics of the system depends on many concrete implementation 
details, \EG, how the situations when the power flow problem has several 
solutions or does not have a solution satisfying all load demands are 
resolved. The dynamics may be altered by the presence of various control 
devices. We will use the equation
 \BEA
   && \DDOT\PHI + D \DOT\PHI = P_{\rm mech} - {\cal P}(P\s \PHI) , 
\label{eq_phi_calP} \\
   && {\cal P}(P\s \PHI) = \underbrace{\frac{P_{\rm supp}}{2}}_{{\sf C} 
x} + \underbrace{\frac{{\sf S}}2 \Bigg( {\sf C} + \sqrt{{\sf C}^2 - 
\frac{3 P_{\rm supp}}{2} - \frac{P_{\rm supp}^2}{{\sf C}^2}} 
\Bigg)}_{{\sf S} y} , \NN
 \EEA where $P_{\rm supp} = {\rm min} \big( P\s {\sf C}^2 / 2 \big)$ --- 
the power demand greater than ${\sf C}^2 / 2$ cannot be fulfilled. The 
``$+$'' sign at the square root in the expression for ${\cal P}$ is, let 
us say, a reflection of some kind of voltage control, where higher 
voltage (larger $y$) is chosen from the two possible solutions.

\begin{figure}
  \begin{picture}(252,140)(0,0)
    \put(0,0){\includegraphics[width=252pt]{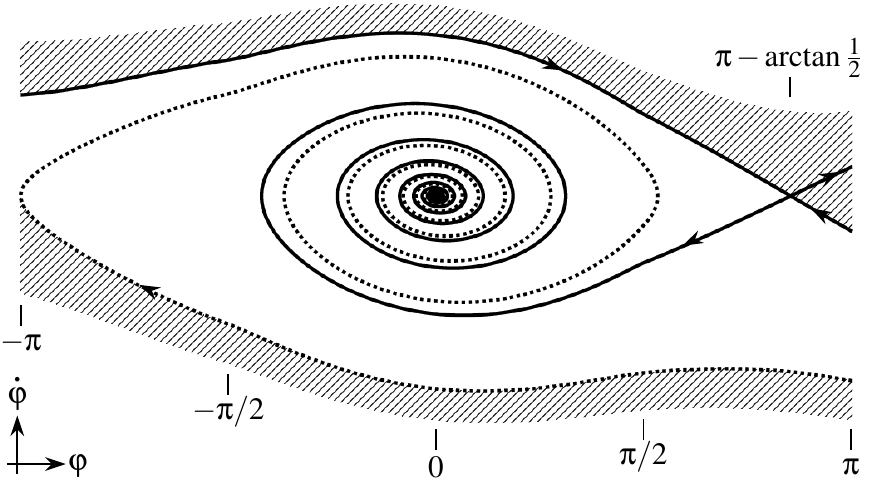}}
  \end{picture}
  \caption{Phase portrait of the system \eqref{eq_phi_calP} with $P(t) 
\equiv P_* = 1 / 4$. Here and in all other $(\PHI\s \DOT\PHI)$-plots the 
scales of $\PHI$ and $\DOT\PHI$ are the same and $D = 0.1$. Shading 
shows the regions where the dynamics eventually leads to $\PHI = \pi$ or 
$\PHI = -\pi$. \label{phase_p14}} \end{figure}

The phase portrait of the system for the case of power demand $P(t) 
\equiv P_* = 1 / 4$ is shown in Fig.~\ref{phase_p14}. The NOP $\PHI = 
\DOT\PHI = 0$ is a stable focus. When $|\PHI| > \pi / 2$ we have $P_{\rm 
supp} = {\sf C}^2 / 2$ and $V = y - \i \, x = {\sf C} \, (2 - \i) / 4$. 
The voltage at the load is small when $\PHI$ is close to $\pi$ or 
$-\pi$.

\begin{figure}
  \begin{picture}(252,124)(0,0)
    \put(0,0){\includegraphics[width=252pt]{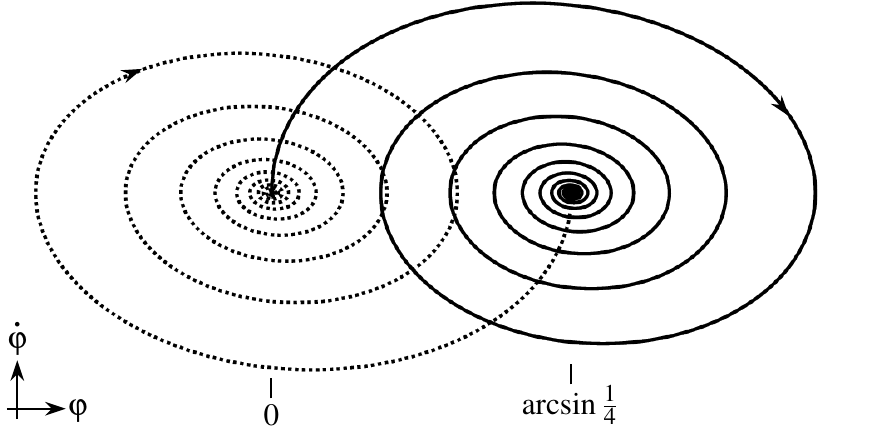}}
  \end{picture}
  \caption{Response of the system to a sudden change of consumed power: 
trajectories with $P = 0$ starting from $\PHI = \DOT\PHI = 0$ (solid 
curve) and with $P = P_* = 1 / 4$ starting from $\PHI = \arcsin 
\frac14$, $\DOT\PHI = 0$ (dotted curve). \label{phase_p0}} \end{figure}

The system is stable to sudden but permanent changes of the power 
demand: \EEG, if $P$ is suddenly dropped to $0$, then the dynamics just 
goes to the new fixed point, see Fig.~\ref{phase_p0}.

\section{Optimal fluctuation equations} \label{sec:equations}

Consider the case when the power demand $P(t)$ is a random process or 
can be changed in time in a controlled way. How the power system can be 
brought out of the stable equilibrium $\PHI = \DOT\PHI = 0$, $P = P_*$ 
to a voltage collapse $\PHI = \pm \pi$? (Even smaller angles $\PHI$ 
could be viewed as an unacceptable loss of synchronization. In this 
work, as an example, will consider the model system ``failure'' being 
$\PHI = \pm \pi$.)

\begin{figure}
  \begin{picture}(252,138)(0,0)
    \put(0,0){\includegraphics[width=252pt]{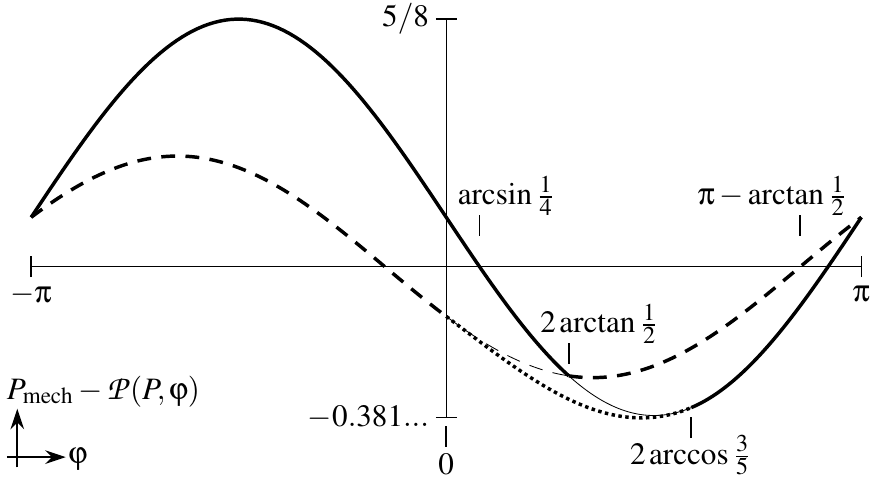}}
  \end{picture}
  \caption{Extreme possible values of $P_{\rm mech} - {\cal P}(P, 
\PHI)$, that are realized at $P = 0$ (solid curve), $P \ge {\sf C}^2 / 
2$ (dashed curve), and $P = {\sf C}^2 (5 {\sf C} - 3) / 4$ (dotted 
curve). \label{value_pin}} \end{figure}

\begin{figure}
  \begin{picture}(252,124)(0,0)
    \put(0,0){\includegraphics[width=252pt]{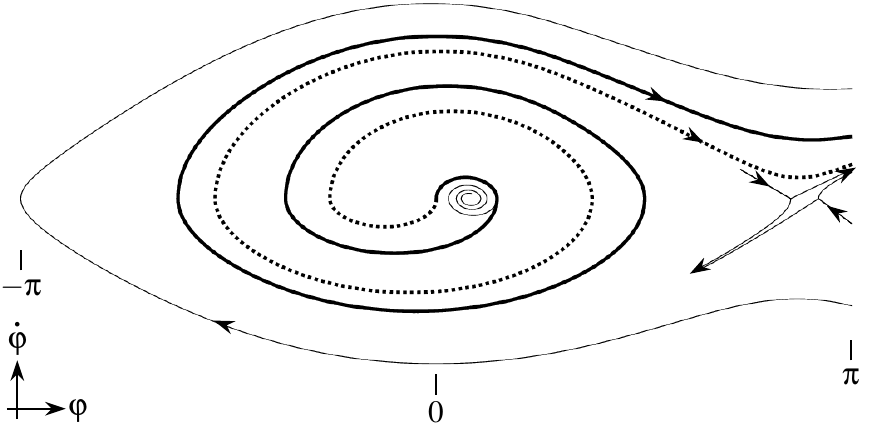}}
  \end{picture}
  \caption{Phase portrait of the system \eqref{eq_phi_calP} with $P$ at 
each moment of time being chosen in order to have maximal if $\DOT\PHI > 
0$ and minimal if $\DOT\PHI < 0$ possible value of thrust $P_{\rm mech} 
- {\cal P}(P, \PHI)$ (see Fig.~\ref{value_pin}). Thick solid{\SEP}dotted 
curve corresponds to starting condition being NOP $\PHI = 0$, $\DOT\PHI 
= 0$, with the initial thrust in positive{\SEP}negative direction. Thin 
spiral near the center is the [solid curve] trajectory from 
Fig.~\ref{phase_p0} shown here for comparison. \label{phase_pin}} 
\end{figure}

Somewhat barbaric but reliable way is to choose $P(t)$ with 
minimal{\SEP}maximal value of ${\cal P}(P\s \PHI)$ when $\DOT\PHI > 
0${\SEP}$\DOT\PHI < 0$. How in this case the power demand $P$ should be 
chosen is shown in Fig.~\ref{value_pin}. Besides the case of $\DOT\PHI < 
0$, $0 < \PHI < 2 \arccos \frac{3}{5} = 4 \arctan \frac{1}{2}$, the 
demand $P$ is equal to extreme supported values: $P = 0$ or $P = {\sf 
C}^2 / 2$. The resulting trajectories leading to $\PHI = \pi$ are shown 
as thick curves in Fig.~\ref{phase_pin}.

The Most Probable Failure Scenario (MPFS) depends on the statistics of 
the random process $P(t)$. (When using optimal fluctuation method for 
practical applications this should be kept in mind, and as much 
properties of the system's random component should be extracted 
experimentally as seems feasible.\footnote{The statistics of the state 
of a power grid at a given time can be written through a functional 
integral \cite{1973_PRA_MSR}. This is especially useful when a random 
noise in the system is correlated in time. When finding the [small] 
probability of failure, the integral often can be estimated by the 
method of steepest descent, with MPFS trajectory being the saddle 
point.}) In this work, as an example, we will use probably the most 
simple yet not entirely unrealistic probabilistic model for $P(t)$: 
There are $\Lambda$ customers, and during each small chunk of time 
$\delta$ each customer either consumes electricity or not, independently 
of other customers and other moments of time. For each chunk of time the 
power demand is binomially distributed. We will assume that $\Lambda \gg 
\lambda \gg 1$, where $\lambda$ is the average number of customers 
consuming electricity at a given time. The binomial distribution is then 
approximated by Poisson: ${\sf P}(n) = \lambda^n \, \e^{-\lambda} / n! 
\approx\footnote{Stirling's approximation: if $n$ is large, then $n! 
\approx \sqrt{2 \pi n} \, (n / \e)^n$.} \exp \big( -\lambda H(n / 
\lambda) \big) / \sqrt{2 \pi n}$ with convex Cram\'er function $H(x) = 1 
- x + x \, \ln x$.

If we track the evolution of the system over time interval $T$, then 
there are $T / \delta$ chunks of time. The probability that there are 
$n_k$ active customers in $k^{\rm th}$ chuck of time is
 \begin{align}
   \prod_{k = 1}^{T / \delta} {\sf P}(n_k) \approx \exp \bigg( -\lambda 
\sum_{k = 1}^{T / \delta} H(n_k / \lambda) \bigg) \bigg/ \prod_{k = 
1}^{T / \delta} \sqrt{2 \pi n_k} . \label{probability}
 \end{align} When $\delta$ is small, the sum over $k$ in the exponent 
could be approximated by an integral over time.

The problem of finding out which failure development is the most 
probable boils down to the following problem of optimal control: Find 
the power consumption $P(t)$, $t_0 \le t \le t_1$ such that $\PHI(t_0) = 
\omega(t_0) = 0$ (NOP) and $|\PHI(t_1)| = \pi$ (failure), and the 
integral ${\cal S} = \int_{t_0}^{t_1} \d t \; H \big( P(t) / P_* \big)$ 
has the minimal possible value. We can fix the time of failure $t_1$ to 
$0$. As we do not care how fast the transition from NOP to the failure 
happens, we have to consider $t_0 = -\infty$.

Consider an equation of motion $\DOT{\VEC x} = {\VEC f}({\VEC x}\s {\VEC 
u})$. We want to move from ${\VEC x}_{\rm start}$ to ${\VEC x}_{\rm 
end}$ in such a way that the integral over time of $H({\VEC u})$ is 
minimal. Let us construct the functional
 \BEA
   {\cal S}\{ \!\!\overbrace{{\VEC x}\s {\VEC u}\s \alpha\s {\VEC 
\beta}}^{{\rm functions~of~}\tau}\!\!\s {\VEC \gamma}_{0, 1} \}
     = {\VEC \gamma}_0 \big( {\VEC x}(\tau_0) - {\VEC x}_{\rm start} 
\big)
     + {\VEC \gamma}_1 \big( {\VEC x}(\tau_1) - {\VEC x}_{\rm end} \big)
   \phantom{\rule{12pt}{1pt}} \hskip-24pt &&
 \NN \\
   + \int\limits_{\tau_0}^{\tau_1} \d \tau \; 
{\VEC \beta}(\tau) \Big( \frac{\d {\VEC x}}{\d \tau} - \alpha {\VEC 
f}({\VEC x}\s {\VEC u}) \Big)
     + \int\limits_{\tau_0}^{\tau_1} \d \tau \; \alpha(\tau) H 
\big( {\VEC u}(\tau) \big) , \hskip-24pt && \NN
 \EEA where ${\VEC \beta}(\tau)$, ${\VEC \gamma}_{0, 1}$ are Lagrange 
multipliers enforcing the equation for ${\VEC x}$ and boundary 
conditions, respectively. The function $\alpha(\tau)$ is introduced to 
arbitrarily parametrize the time $t$ ($\d t = \alpha \, \d \tau$), this 
trick\footnote{We could introduce ${\cal S}$ with integral over time $t$ 
from $0$ to $T$ (or make $\alpha \equiv 1$, $\tau = t$, and $\tau_0 = 
0$, $\tau_1 = T$), and then try to find the minimum in $T$. This may be 
technically cumbersome as it would involve derivatives of ${\VEC x}(t)$ 
and ${\VEC u}(t)$ with respect to $T$.} allows us to have fixed limits 
of integration $\tau_{0, 1}$ (we can choose any values, \EG, $\tau_0 = 
0$ and $\tau_1 = 1$).

Setting the variation over various variables in in ${\cal S}$ to zero we 
get Euler--Lagrange equations
 \begin{align}
   {\VEC\beta} & ~~:~~ \DOT{\VEC x} = {\VEC f}({\VEC x}\s {\VEC u}) , 
\label{eq_beta} \\
   {\VEC x} & ~~:~~ \DOT\beta_i = -{\VEC\beta} \cdot {\partial {\VEC f}} 
/ {\partial x_i} , \label{eq_x} \\
   \alpha & ~~:~~ H({\VEC u}) - {\VEC\beta} \cdot {\VEC f}({\VEC x}\s 
{\VEC u}) = 0 , \label{eq_alpha} \\
   {\VEC u} & ~~:~~ {\VEC u}\mbox{ minimizes }H({\VEC u}) - {\VEC\beta} 
\cdot {\VEC f}({\VEC x}\s {\VEC u}) . \label{eq_u}
 \end{align} From the structure of the expression for ${\cal S}$ another 
view on $\alpha$ is that it is a Lagrange multiplier enforcing 
eq.~\eqref{eq_alpha}.

If one considers the dynamics in $({\VEC x}\s {\VEC\beta}\s {\VEC u})$ 
space (of dimension $2 D_x + D_u$, where $D_x = D_\beta$ and $D_u$ are 
the lengths of vectors ${\VEC x}$ and ${\VEC u}$; in our model system 
$D_x = 2$ and $D_u = 1$), then the eqs.~\eqref{eq_alpha} and 
\eqref{eq_u} produce a $(2 D_x - 1)$-dimensional surface in it. The 
equations of motion \eqref{eq_beta} and \eqref{eq_x} determine the 
direction of movement, which gives another $(2 D_x - 1)$ conditions. 
From naive dimensions counting one could think that together all the 
equations are satisfied at isolated points. This is not true.

In a ``good'' case, when all functions are differentiable, because of 
the eqs.~\eqref{eq_beta}, \eqref{eq_x}, and \eqref{eq_u} the quantity 
${\cal H} = H({\VEC u}) - {\VEC\beta} \cdot {\VEC f}({\VEC x}\s {\VEC 
u})$ is conserved, $\d {\cal H} / \d t = 0$, so the eq.~\eqref{eq_alpha} 
is automatically kept true, and all the equations are satisfied along 
$1$-dimensional [optimal] trajectories. The eq.~\eqref{eq_u} is 
Pontryagin's minimum principle with ${\cal H}$ being the 
Hamiltonian\footnote{We have ${\DOT x}_i = -\partial {\cal H} / \partial 
\beta_i$ and $\DOT\beta_i = \partial {\cal H} / \partial x_i$.}, the 
eq.~\eqref{eq_x} is the costate equation, and the eq.~\eqref{eq_alpha} 
reflects that the evolution time is not fixed \cite{Pontryagin}.

For our model power grid we have
 \begin{align}
   \frac{\d}{\d t} \underbrace{\left[ \begin{array}{c} \PHI \\ \omega 
\end{array} \right]}_{\VEC x} &= \underbrace{\left[ \begin{array}{c} 
\omega \\ -D \omega + P_{\rm mech} - {\cal P}(P\s \PHI) \end{array} 
\right]}_{{\VEC f}({\VEC x}\s {\VEC u})} , \tag{\ref{eq_beta}a} 
\label{eq_beta_system} \\
   \frac{\d}{\d t} \underbrace{\left[ \begin{array}{c} \beta_\PHI \\ 
\beta_\omega \end{array} \right]}_{\VEC\beta} &= \left[ \begin{array}{c} 
\beta_\omega \; \partial {\cal P}(P\s \PHI) / \partial \PHI \\ D 
\beta_\omega - \beta_\PHI \end{array} \right] , \tag{\ref{eq_x}a} 
\label{eq_x_system} \\
  \beta_\PHI = \frac{1}{\omega} \Big( H( &P / P_*) + \beta_\omega \big( 
D \omega - P_{\rm mech} + {\cal P}(P\s \PHI) \big) \Big) . 
\tag{\ref{eq_alpha}a} \label{eq_alpha_system}
 \end{align} Let us discuss two issues with these equations. 
\emph{First}, sometimes the eq.~\eqref{eq_u} prescribes $P = {\sf C}^2 / 
2$, and then the derivative $\partial {\cal P}(P\s \PHI) / \partial 
\PHI$ does not exist (or is infinite). In this case $\beta_\PHI$ could 
rapidly change in time, and its value should be chosen so that the 
eq.~\eqref{eq_alpha_system} is maintained.

\emph{Second}, if $\omega = 0$, then how the eq.~\eqref{eq_alpha_system} 
should be understood? When $\omega = 0$ the movement in the $(\PHI\s 
\DOT\PHI)$ phase plane is directed vertically no matter what is the 
value of $P$. As the time of the evolution is not important for us, we 
don't care how fast we are moving if the direction of movement is the 
same. Thus, in order to minimize $\int \d t \; H(P / P_*)$ we have to 
choose $P = P_*$. On the MPFS trajectory we have $P(t) = P_* \big( 1 + 
\omega(t) \rho(t) \big)$ with $\rho$ being finite at $\omega = 0$. To 
get that from the eq.~\eqref{eq_u} we need to have $\beta_\omega = 
\omega \beta$ with finite $\beta$. The eq.~\eqref{eq_alpha_system} is 
then non-singular\footnote{The change of variables $P$, $\beta_\omega$ 
$\to$ $\rho$, $\beta$ made by division by $\omega$ is an example of 
singularity resolution by ``blowing up''.} at $\omega = 0$. The MPFS 
equations become
 \begin{align}
   \DOT\PHI &= \omega , \quad
   \DOT\omega = -D \omega + P_{\rm mech} - {\cal P} \big( P_* (1 + 
\omega \rho)\s \PHI \big) , \tag{\ref{eq_beta}b} \label{eq_beta_simple} 
\\
   \DOT\beta &= D \beta - H(1 + \omega \rho) / \omega^2 , 
\tag{\ref{eq_x},\ref{eq_alpha}b} \label{eq_x_simple} \\
   \rho&\mbox{ minimizes }H(1 + \omega \rho) + \omega \beta {\cal P} 
\big( P_* (1 + \omega \rho)\s \PHI \big) . \tag{\ref{eq_u}b} 
\label{eq_u_simple}
 \end{align} When $\beta$ is small, then due to the 
eq.~\eqref{eq_u_simple} $\rho$ is small too, and in the 
eq.~\eqref{eq_x_simple} the value of $\beta$ never changes sign. In MPFS 
$\beta$ is positive. Positive{\SEP}negative $\omega \beta$ (or just 
$\omega$) bias the thrust $P_{\rm mech} - {\cal P}(P\s \PHI)$ to 
larger{\SEP}smaller values, and that causes NOP to be unstable. When 
$|\omega \beta| \gg 1$ the extreme values of thrust, as in 
Fig.~\ref{value_pin}, are realized.

\section{Numerical solution from IVP} \label{sec:IVP}

The MPFS equations (\ref{eq_beta_simple}--\ref{eq_u_simple}) form a 
non-linear Boundary Value Problem (BVP). Here we show how it can be 
treated as an Initial Value Problem (IVP). (This can be viewed as a kind 
of analog of the shooting (from $t = -\infty$) method.) Let us consider 
the linearization of MPFS equations near the NOP solution $\PHI = \omega 
= 0$:
 \begin{align}
   \DOT\PHI &= \omega , \quad
   \DOT\omega = -D \omega - \omega \rho / 8 - 7 \PHI / 16 , 
\tag{\ref{eq_beta}c} \label{eq_beta_lin} \\
   \DOT\beta &= D \beta - \rho^2 / 2 , \quad
   \rho = -\beta / 8 . \tag{\ref{eq_x},\ref{eq_alpha},\ref{eq_u}c}
 \end{align} There are two steady state solutions: $\beta = \rho = 0$, 
\IE, $P = P_*$ (stable in $\PHI$ and $\omega$, unstable in $\beta$ 
direction, this solution is not interesting for us\footnote{The 
equations (\ref{eq_beta}--\ref{eq_u}) for $\DOT{\VEC x} = {\VEC f}({\VEC 
x}\s {\VEC u})$ support solution ${\VEC\beta} \equiv {\VEC 0}$, ${\VEC 
u}_*$ minimizes $H({\VEC u})$ with $H({\VEC u}_*) = 0$. It works if we 
reach ${\VEC x}_{\rm end}$ from ${\VEC x}_{\rm start}$ in a relaxed 
manner, \IE, with ${\VEC u} \equiv {\VEC u}_*$.}); and $\beta = 128 D$, 
$\rho = -16 D$ which will give us the MPFS trajectory.

\begin{figure}
  \begin{picture}(252,115)(0,0)
    \put(0,0){\includegraphics[width=252pt]{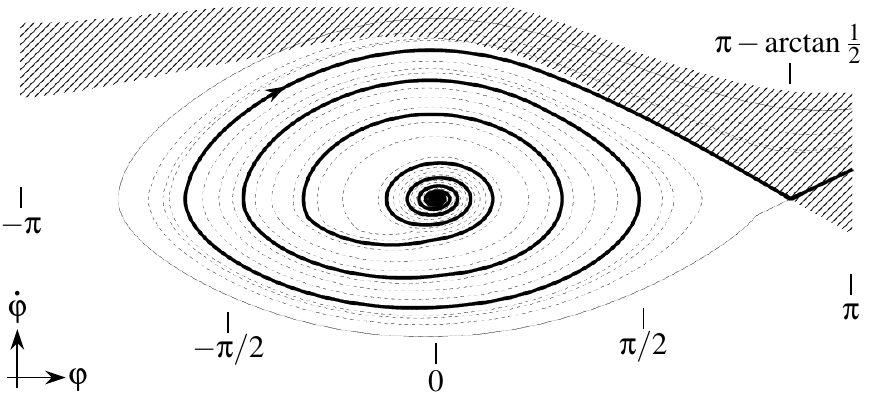}}
  \end{picture}
  \caption{MPFS trajectory constructed with $\EPS = \EPS_0 \approx 
1.03452 \cdot 10^{-12}$. Three dashed thin line trajectories are built 
with $\EPS = \EPS_0 F^{1 / 4}$, $\EPS_0 F^{1 / 2}$, and $\EPS_0 F^{3 / 
4}$, where $F = \exp \big( 2 \pi D / \sqrt{7 / 4 - D^2} \big)$ 
(trajectory with $\EPS = \EPS_0 F$ would almost coincide with the one 
with $\EPS = \EPS_0$). The thin line trajectory coming from the saddle 
at $\PHI = \pi - \arctan \frac12$ uses $\EPS$ that is slightly less than 
$\EPS_0$. \label{ivp_instanton}} \end{figure}

\begin{figure}
  \begin{picture}(252,92)(0,0)
    \put(0,0){\includegraphics[width=252pt]{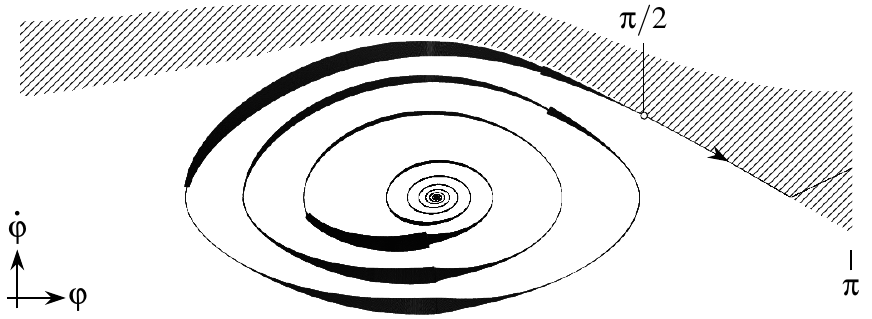}}
  \end{picture}
  \caption{MPFS trajectory constructed with $\EPS = \EPS_0 \approx 
1.03452 \cdot 10^{-12}$. The thickness of the line is equal to $\big( 
H(P / P_*) + 0.01 \big) \cdot 12\mbox{ pt}$. The trajectory point at 
$\PHI = \pi / 2$ is right on the border of the shaded region. 
\label{thickness}} \end{figure}

\begin{figure}
  \begin{picture}(252,98)(0,0)
    \put(0,0){\includegraphics[width=252pt]{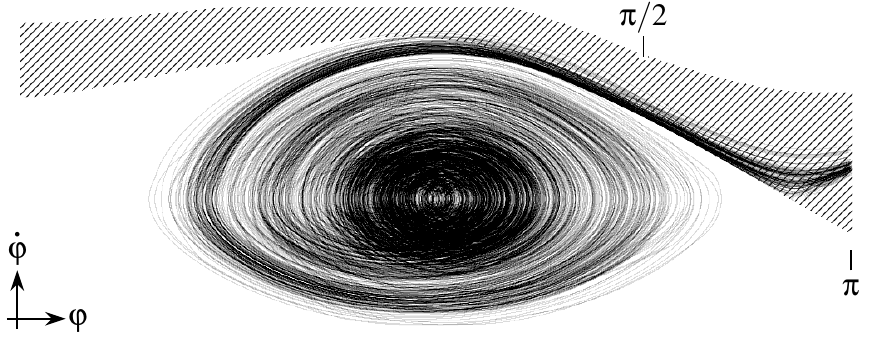}}
  \end{picture}
  \caption{$100$ failure trajectories obtained by direct numerical 
simulation of \eqref{eq_phi_calP} with $\lambda = 2$ and $\delta = 1.5$ 
(mean time to failure is about $4.8 \cdot 10^7$). On each trajectory the 
last $100$ time units before reaching $\PHI = \pi$ are shown. 
\label{statistics}} \end{figure}

\begin{figure}
  \begin{picture}(252,130)(0,0)
    \put(0,0){\includegraphics[width=252pt]{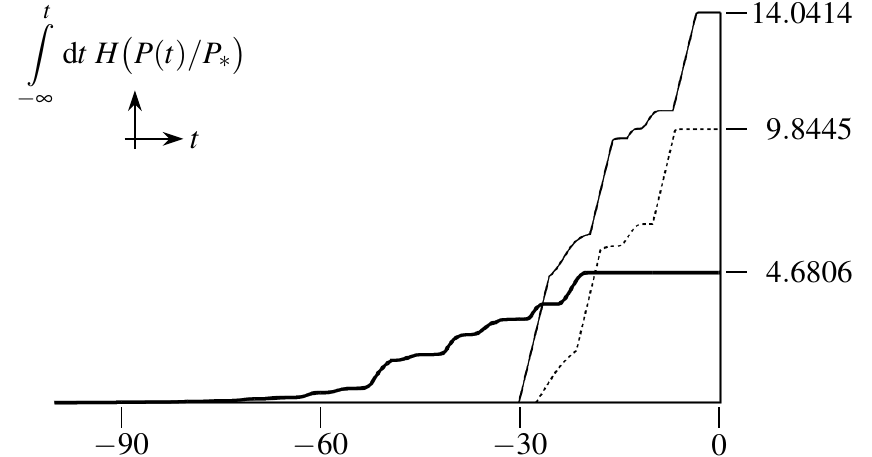}}
  \end{picture}
  \caption{The growth of cost function $S(t) = \int_{-\infty}^t \d t \; 
H \big( P(t) / P_* \big)$ in time. The value of ${\cal S} = S(0)$ is 
indicated by numbers on the right. Thin solid and dashed curves 
correspond to ``barbaric'' trajectories from Fig.~\ref{phase_pin}. The 
plateau on the thick curve near $t = 0$ is technically infinitely wide, 
as MPFS trajectory passes through the saddle at $\PHI = \pi - \arctan 
\frac12$, $\DOT\PHI = 0$. \label{action}} \end{figure}

We solve (until failure) the IVP 
eqs.~(\ref{eq_beta_simple}--\ref{eq_u_simple}) with the initial data 
$\PHI_{\rm start} = 0$, $\omega_{\rm start} = \EPS \ll 1$, and 
$\beta_{\rm start} = 128 D$. Initially, while $\rho \approx -16 D$, the 
eqs.~\eqref{eq_beta_lin} look like $\DOT\PHI = \omega$, $\DOT\omega = D 
\omega - 7 \PHI / 16$, and NOP $\PHI = \omega = 0$ is an unstable focus 
--- the IVP solution is periodic in $\log \EPS$ (this can be interpreted 
as going from one swirl to the next one in the spiral). We need to cover 
just one period to get all solutions. From this $1$-dimensional family 
of solutions we need to pick the one with the minimal value of $\int \d 
t \; H \big( P(t) / P_* \big)$.

The resulted MPFS trajectory is shown in Figs.~\ref{ivp_instanton} and 
\ref{thickness}, see also Fig.~\ref{statistics}. In Fig.~\ref{action} 
one can see the buildup of ${\cal S} = \int_{-\infty}^0 \d t \; H(P / 
P_*)$ in time and comparison of ${\cal S}$ values on the optimal 
trajectory and the ones from Fig.~\ref{phase_pin}.

Obtaining the solution from IVP is not going to work well for large 
systems. When the dimension $D_x$ of vector ${\VEC x}$ is large, it is 
hard to pick the most probable trajectory from $(D_x - 1)$-dimensional 
family of IVP solutions.

\section{Hamilton--Jacobi--Bellman equation} \label{sec:HJB}

Let us introduce the following value function: ${\cal S}(\PHI\s 
\DOT\PHI;\, T)$ is the minimal value of integral of $H \big( P(t) / P_* 
\big)$ over time on the trajectory that starts at $(\PHI\s \DOT\PHI)$ 
and reaches $\PHI = \pm\pi$ during the time {\emph{less or equal}} to 
$T$. (If such a trajectory does not exist, then let us define ${\cal S} 
= +\infty$.) This function satisfies the Hamilton--Jacobi--Bellman (HJB) 
equation.

We are interested in ${\cal S}(\PHI\s \DOT\PHI) = \lim_{T \to \infty} 
{\cal S}(\PHI\s \DOT\PHI;\, T)$ (we care about the overall probability 
of the failure, not in how much time is needed for it to develop); and 
more specifically in ${\cal S}(0\s 0)$, which corresponds to NOP as the 
starting state. The \emph{time-independent} HJB equation in this case 
are the eqs.~\eqref{eq_alpha} and \eqref{eq_u}, with $-{\VEC\beta}$ 
being the gradient of the value function ${\cal S}(\PHI\s \DOT\PHI)$. 
The eq.~\eqref{eq_alpha} determines the rate of change of $\int \d t \; 
H(P / P_*)$ through this gradient and the equation of motion.

From Fig.~\ref{value_pin} we see that when $\PHI > 2 \arctan \frac12$ 
the most effective way to push the system to the right is to choose $P 
\ge {\sf C}^2 / 2$. If $\PHI > \pi / 2$, then ${\sf C}^2 / 2 < P_*$, so 
$P = P_*$ is both the most effective and the most probable. If the 
system didn't reach the upper shaded region in Fig.~\ref{phase_p14} 
before $\PHI = \pi / 2$, then it will not be able to reach the saddle at 
$\PHI = \pi - \arctan \frac12$ and will be pushed back to the left. The 
optimal trajectory should pass through the border point of the shaded 
region at $\PHI = \pi / 2$, see also Fig.~\ref{thickness}.

The function ${\cal S}(\PHI\s \DOT\PHI)$ is not continuous and undergoes 
a jump along the $\pi / 2 < \PHI < \pi$ part of the border of the shaded 
region. The height of the jump is the integral of $H(P / P_*)$ along the 
additional [penalty] loop. In power grids with bounds on the possible 
dynamics due to, \EG, existence of power flow solution regions or 
maximum allowed currents through lines, such discontinuities should be a 
commonplace.

In general situation a time-independent value function can be introduced 
if the random process describing fluctuations in the system (\EG, 
fluctuations of power demand or generation in renewable energy sources) 
is Markovian (in the above example the values of the demand at different 
times are independent) --- it is then the function of both the system's 
state and a state of the Markov random process.

\section{Numerical solution from BVP} \label{sec:BVP}

There are two major ways to solve optimal control problems --- 
[historically earlier] indirect \cite{Pontryagin} and direct 
\cite{Betts} methods. Here we show how to obtain the MPFS trajectory by 
indirect method, \IE, by solving the BVP 
(\ref{eq_beta_simple}--\ref{eq_u_simple}). There are some similarities 
with multiple shooting and Gauss--Seidel methods.

If one tries to solve the BVP for the whole trajectory by 
quasi-linearization method or [multiple-]shooting with Newton's method, 
then the initial guess should be too close the the solution, otherwise 
Newton's iterations do not converge. Continuation method, with slow 
dragging of the final point to $\PHI = \pi$, helps very little --- the 
drag should be quite slow. The possible reason is that the whole MPFS 
trajectory has rich enough history with many features.

A way to overcome these difficulties is to solve BVP for small parts of 
the trajectory, such BVPs should be much simpler to deal with. The end 
points of the trajectory segment should be such that it is possible to 
go from ${\VEC x}_{\rm start}$ to ${\VEC x}_{\rm end}$ in a 
straightforward manner for some realization of ${\VEC u}(t)$. This is 
achieved by generating a suitable (\IE, leading to failure) seed 
trajectory with a certain ${\VEC u}(t)$, which is then gradually 
improved by optimizing its [not large] segments.

\begin{figure}
  \begin{picture}(252,115)(0,0)
    \put(0,0){\includegraphics[width=252pt]{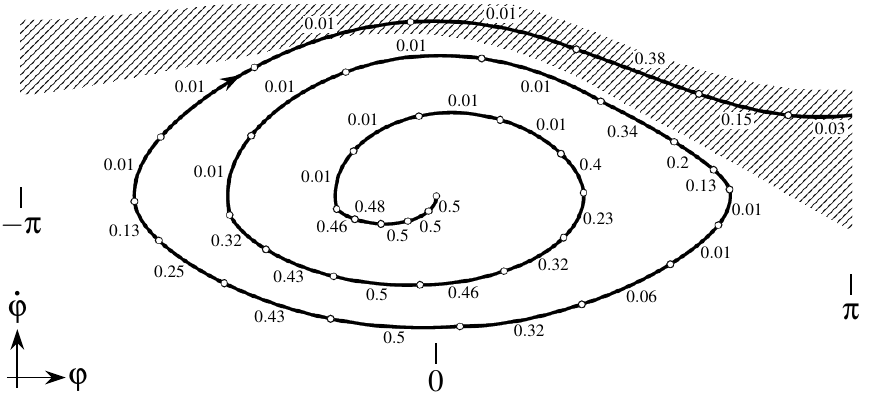}}
  \end{picture}
  \caption{Seed trajectory made of segments with duration $1$. The value 
of the power demand $P$ on each segment is written next to it.
 \label{seed}} 
\end{figure}

We generate the seed trajectory starting from NOP $\PHI = \DOT\PHI = 0$ 
and adding to it one by one segments with constant value of $P$ and with 
duration $1$. If one chooses the new segment's value of $P$ such that 
the value of $\PHI_{\rm end}$ at the right end of the segment is 
maximal, then one gets the solid curve trajectory from 
Fig.~\ref{phase_p0}. The choice of the value of $P$ such that the change 
of $\PHI$ over the segment, $|\PHI_{\rm end} - \PHI_{\rm start}|$, is 
maximal gives the seed trajectory\footnote{In general situation the 
creation of seed trajectory could be tricky, and additions of several 
new segments at once (checking all possible combinations of segments' 
control values) could be needed.} shown in Fig.~\ref{seed}. Initially we 
set $\beta(t) \equiv 0$.

At the very start of the trajectory we have $\PHI = \omega = 0$. In the 
system of eqs.~(\ref{eq_beta_simple}--\ref{eq_u_simple}) no matter what 
is the starting value of $\beta$, we have $P = P_*$ and don't leave NOP. 
We ``jump-start'' the trajectory by adding to it the solution of the BVP 
eqs.~(\ref{eq_beta_system},\ref{eq_x_system},\ref{eq_u}), that is better 
defined near NOP. Its duration should be large enough (we used $50$). 
The boundary conditions are that we start at NOP and end up at some 
point of the seed trajectory that is close to NOP (we chose a point with 
$0.2$ distance from NOP) --- the power demand $P$ is then far from ${\sf 
C}^2 / 2$, and the eq.~\eqref{eq_alpha} is kept true.

In accordance with Bellman's principle of optimality \cite{Bellman} at 
the end of the trajectory we should have $P = P_*$ --- if $\DOT\PHI > 0$ 
we'll reach $\PHI = \pi$ anyway. We ``relax'' the trajectory by checking 
from where we can set $P = P_*$ and still get a failure.

Next we repeatedly choose a segment within the trajectory and attempt to 
optimize it. The segment's BVP, 
eqs.~(\ref{eq_beta_simple}--\ref{eq_u_simple}), is solved by simple 
shooting with Newton's method. The boundary conditions fix $\PHI_{\rm 
start}$, $\omega_{\rm start}$, $\PHI_{\rm end}$, and $\omega_{\rm end}$ 
--- they are not changed. The shooting parameters are $\beta_{\rm 
start}$ and the new duration of the segment. It is not crucial that the 
segment's BVP is solved right at the moment, so the maximal number of 
iterations $N_{\rm iter}$ in Newton's method is not needed to be large 
(we used $N_{\rm iter} = 10$). If the solution of the BVP is found, then 
it substitutes the old state of the segment. As the number of grid 
points could easily change, a doubly-linked list as a data structure for 
the trajectory is convenient.

The shorter is the segment, the smaller is its change due to 
optimization. If we work only with short segments, then the progress 
will be very slow. Thus, if obtaining the solution for the segment's BVP 
was successful, we try the segment with the same starting point (with 
already updated value of $\beta_{\rm start}$) but longer (we multiplied 
the duration of the segment by $1.25$). Sometimes even if solving the 
BVP for a very short segment fails, it is beneficial to try a long 
segment anyway --- it could throw a bridge over difficult trajectory 
parts (like the first $5$ segments of the seed trajectory on 
Fig.~\ref{seed}, where $P \ge {\sf C}^2 / 2$ results in the extreme 
values of $P_{\rm mech} - {\cal P}(P, \PHI)$ which demands values of 
$\beta$ to be large\footnote{Such artifacts in seed trajectory should be 
avoided if possible.}).

If we fail to obtain the solution, then we shift the segment's starting 
point forward by a share of its length (we shifted to a random point in 
the first half of the segment) and try again. The value of $\beta$ could 
be already non-zero [and close to being reasonable] from previous 
successes. Eventually we complete a ``pass'' --- go through the whole 
trajectory. We continue to run passes until the trajectory converges.

\begin{figure}
  \begin{picture}(252,204)(0,0)
    \put(0,0){\includegraphics[width=252pt]{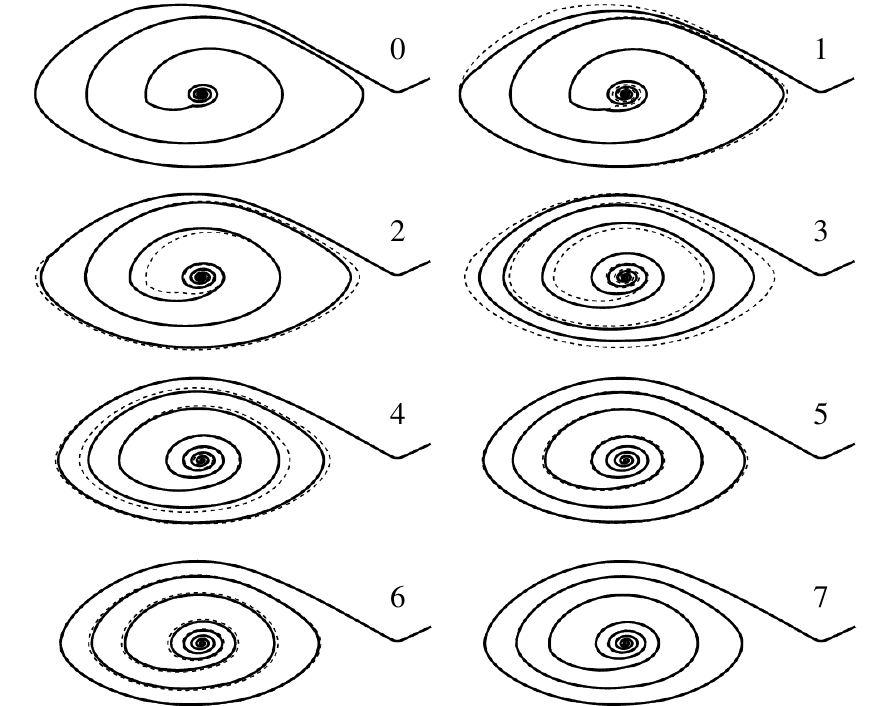}}
  \end{picture}
  \caption{The progress in obtaining the BVP solution. Attached to each 
trajectory is the number of completed passes. In upper left is the 
``jump-started'' and ``relaxed'' seed trajectory from Fig.~\ref{seed}. 
The dashed thin curve provides the comparison with the trajectory from 
the previous pass. \label{bvp_instanton}} \end{figure}

The gradual convergence of the trajectory to the optimal solution is 
shown in Fig.~\ref{bvp_instanton}.

\section{Probability of failure} \label{sec:probability}

When both $\PHI$ and $\omega$ are small, and $P \approx P_*$, we can 
write the linearized equations of motion near NOP as
 \begin{align}
   \DDOT\PHI + D \DOT\PHI = - (P - P_*) / 2 - 7 \PHI / 16 . 
\label{stochastic_lin}
 \end{align} When the average number of consuming customers is large, 
$\lambda \gg 1$, and the with of the time chunk is small, $\delta \ll 
1$, we can think about $P - P_*$ as a noise with small correlation time. 
It has (as ${\sf E} \, P(t) = P_*$, ${\sf E} \, P^2(t) = P_*^2 (1 + 1 / 
\lambda)$) zero mean and the integral of its autocovariance function 
over time is equal to $P_*^2 \delta / \lambda$. This results in 
probability distribution of the system's state being mainly a Gaussian 
distribution with NOP $\PHI = \DOT\PHI = 0$ as its center. After some 
calculations we get the probability density function $f(\PHI\s \DOT\PHI) 
\propto \exp \big( -(\lambda / \delta) \, 4 D (7 \PHI^2 + 16 \DOT\PHI^2) 
\big)$. This can be obtained, \EG, by directly treating the stochastic 
system \eqref{stochastic_lin} or by solving the linearized MPFS 
equations \eqref{eq_beta_lin} with $\rho = -16 D$ and the boundary 
condition that the trajectory ends at the state $(\PHI\s \DOT\PHI)$.

This found Gaussian center of the probability distribution could be used 
for a naive estimate of the probability of failure. We are interested in 
reaching $\PHI = \pm \pi$, so the probability of failure in unit time is 
${\sf P}(\mbox{failure}) \sim f(\pi\s 0) \sim \exp \big( -(\lambda / 
\delta) \, 28 D \pi^2 \big)$.

From \eqref{probability} we estimate ${\sf P}(\mbox{failure}) \sim \exp 
\big( -(\lambda / \delta) \, {\cal S} \big)$. Of course, $\exp \big( 
-(\lambda / \delta) \, {\cal S} \big)$ is not the realization 
probability of the most probable trajectory that leads to a failure. The 
expression \eqref{probability} contains also the factor $1 / \sqrt{2 \pi 
n}$ for each chunk of time, which makes the probability of that very 
trajectory to be much smaller than $\exp \big( -(\lambda / \delta) \, 
{\cal S} \big)$. But the most probable trajectory is not the only 
trajectory that leads to a failure. Lots of trajectories that are close 
to the MPFS trajectory also cause a failure. To accurately find the 
failure probability, we need to take into account all such trajectories.

Exploring the vicinity of the MPFS trajectory is similar to the study of 
the saddle point in the method of steepest descent. This is typically a 
very hard task. Luckily, we often may skip it. The MPFS trajectory 
(position of the saddle point) provides the exponential part of the 
answer $\exp \big( -(\lambda / \delta) \, {\cal S} \big)$, while taking 
into account the whole bunch of trajectories (integrating in the 
vicinity of the saddle point) provides the pre-exponential correction 
which is less important.

Considering the MPFS trajectory over infinite interval of time 
$(-\infty\s 0)$, we get infinitely many time chunks, and thus $1 / 
\sqrt{2 \pi n}$ factors. At the same time the integral ${\cal S} = 
\int_{-\infty}^0 \d t \; H \big( P(t) / P_* \big)$ converges because the 
integrand exponentially decays when $t$ tends to $-\infty$. The 
probability of the MPFS trajectory realization, strictly speaking, is 
equal to $0$. The probability of failure is not, because in this case 
there are infinitely many different trajectories that end up in failure. 
As $D > 0$, the details of the system evolution in a very distant past 
are not important --- they are forgotten because of the decay. The value 
of $n_k$ of the $k^{\rm th}$ chunk of time in the distant past doesn't 
affect whether the failure did happen or not, and taking into account 
all trajectories that end up in failure means summation over $n_k$ which 
eliminates the factor $1 / \sqrt{2 \pi n_k}$.

The whole part of the MPFS trajectory (which can be seen in 
Figs.~\ref{ivp_instanton} and \ref{thickness}) that is inside the 
[Gaussian] center of the probability distribution is irrelevant --- the 
variation of the trajectory there hardly affects the outcome of the 
evolution. The larger is $\lambda / \delta$, the smaller is this part 
(and the smaller is the probability ${\sf P}(\mbox{failure})$). The 
meaningful part of the MPFS trajectory is the one outside of the region 
where the system spends most of its time.

Although the probability distribution of the system's state has its 
center being close to Gaussian, its tails are not continuing the trend. 
To see how strongly the naive estimation from the Gaussian center is 
wrong, we can compare for $D = 0.1$ the numerical values of ${\cal S} 
\approx 4.6806$ and of $28 D \pi^2 \approx 27.6$. Their ratio is about 
$5.9$, and these numbers are in the exponent. With such a difference, if 
for example we naively estimate the probability of failure in unit time 
as about $10^{-12}$, then the actual probability could be not far from 
$10^{-2}$. This illustrates the danger of estimating the probability of 
rare events from such characteristics of the probability distribution as 
mean value and standard deviation.\footnote{One need not think that 
assumption of Gaussian statistics always underestimates the probability 
of a rare event. Consider, \EG, our model power grid with $D$ being 
large enough, so it is impossible to bring the system to $\PHI = \pm 
\pi$.} Same sentiments against widespread assumption of Gaussian or 
normal statistics were expressed, \EG, in \cite{2004_PD_FKLMT}.

\begin{figure}
  \begin{picture}(252,100)(0,0)
    \put(0,0){\includegraphics[width=252pt]{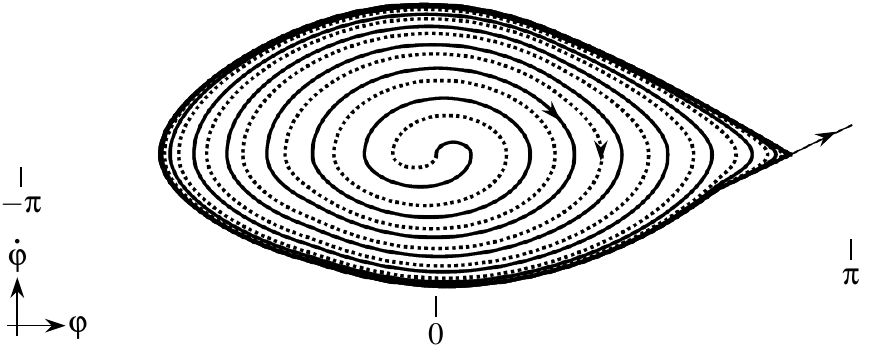}}
  \end{picture}
  \caption{Two trajectories of the system \eqref{eq_phi_calP} with $P$ 
at each moment of time being chosen (with a condition $|P - P_*| < 
P_{\rm control}$) in order to have maximal if $\DOT\PHI > 0$ and minimal 
if $\DOT\PHI < 0$ possible value of thrust $P_{\rm mech} - {\cal P}(P, 
\PHI)$. Thick solid{\SEP}dotted curve corresponds to starting condition 
being NOP $\PHI = 0$, $\DOT\PHI = 0$, with the initial thrust in 
positive{\SEP}negative direction. The used value of $P_{\rm control} 
\approx 0.1385$ is the minimal one for which it is possible to bring the 
system to $\PHI = \pi$ state. \label{phase_pin_control}} \end{figure}

There could be situations where instead of failure probability the 
relevant question is whether a malicious party that controls part of a 
grid could bring the system to failure. For an example, see 
Fig.~\ref{phase_pin_control}.

\section{Conclusions} \label{sec:conclusions}

Even if a power system is ``statically stable'', \IE, no constant in 
time choice of such parameters as power demands at loads or production 
by generators could drive the system to a non-operating state, by 
varying these parameters in time we may be able to bring the system to a 
failure.

While designing a power system, it is interesting to know what is the 
probability per unit time that causing a failure fluctuation would 
occur. This probability can be effectively estimated by optimal 
fluctuation method \cite{1964_AP_L}, in which the most probable 
evolution leading to a failure is a solution to a certain optimal 
control problem.

The equations for the ``optimal'' failure trajectory are sensitive to 
the details of random power demand or generation statistics. For 
accurate estimates of a failure probability, realistic models of various 
power grid components with randomly varying parameters need to be 
developed.

\bibliographystyle{IEEEtran}
\bibliography{refs}

\begin{thebibliography}{10}
\providecommand{\url}[1]{#1}
\csname url@samestyle\endcsname
\providecommand{\newblock}{\relax}
\providecommand{\bibinfo}[2]{#2}
\providecommand{\BIBentrySTDinterwordspacing}{\spaceskip=0pt\relax}
\providecommand{\BIBentryALTinterwordstretchfactor}{4}
\providecommand{\BIBentryALTinterwordspacing}{\spaceskip=\fontdimen2\font plus
\BIBentryALTinterwordstretchfactor\fontdimen3\font minus
  \fontdimen4\font\relax}
\providecommand{\BIBforeignlanguage}[2]{{%
\expandafter\ifx\csname l@#1\endcsname\relax
\typeout{** WARNING: IEEEtran.bst: No hyphenation pattern has been}%
\typeout{** loaded for the language `#1'. Using the pattern for}%
\typeout{** the default language instead.}%
\else
\language=\csname l@#1\endcsname
\fi
#2}}
\providecommand{\BIBdecl}{\relax}
\BIBdecl

\bibitem{1989_SCL_DC}
I.~Dobson and H.-D. Chiang, ``Towards a theory of voltage collapse in electric
  power systems,'' \emph{Syst. Control Lett.}, vol.~13, no.~3, pp. 253--262,
  1989.

\bibitem{2007_C_DCLN}
I.~Dobson, B.~A. Carreras, V.~E. Lynch, and D.~E. Newman, ``Complex systems
  analysis of series of blackouts: Cascading failure, critical points, and
  self-organization,'' \emph{Chaos}, vol.~17, no.~2, p. 026103, 2007.

\bibitem{2009_IEEETPS_SWB}
J.~Salmeron, K.~Wood, and R.~Baldick, ``Worst-case interdiction analysis of
  large-scale electric power grids,'' \emph{IEEE Trans. Power Syst.}, vol.~24,
  no.~1, pp. 96--104, 2009.

\bibitem{1964_AP_L}
I.~M. Lifshitz, ``The energy spectrum of disordered systems,'' \emph{Adv.
  Phys.}, vol.~13, no.~52, pp. 483--536, 1964.

\bibitem{1966_PR_HL}
B.~I. Halperin and M.~Lax, ``Impurity-band tails in the high-density limit.
  {I}. minimum counting methods,'' \emph{Phys. Rev.}, vol. 148, no.~2, pp.
  722--740, 1966.

\bibitem{1966_PR_ZL}
J.~Zittartz and J.~S. Langer, ``Theory of bound states in a random potential,''
  \emph{Phys. Rev.}, vol. 148, no.~2, pp. 741--747, 1966.

\bibitem{1975_PLB_BPST}
A.~A. Belavin, A.~M. Polyakov, A.~S. Schwartz, and Y.~S. Tyupkin,
  ``Pseudoparticle solutions of the {Y}ang--{M}ills equations,'' \emph{Phys.
  Lett. B}, vol.~59, no.~1, pp. 85--87, 1975.

\bibitem{1977_JETP_L}
L.~N. Lipatov, ``Divergence of the perturbation-theory series and the
  quasi-classical theory,'' \emph{Sov. Phys. JETP}, vol.~45, no.~2, pp.
  216--223, 1977.

\bibitem{1996_PRE_FKLM}
G.~Falkovich, I.~Kolokolov, V.~Lebedev, and A.~Migdal, ``Instantons and
  intermittency,'' \emph{Phys. Rev. E}, vol.~54, no.~5, pp. 4896--4907, 1996.

\bibitem{2011_IEEETSG_CPS}
M.~Chertkov, F.~Pan, and M.~G. Stepanov, ``Predicting failures in power grids:
  the case of static overloads,'' \emph{IEEE Trans. Smart Grids}, vol.~2,
  no.~1, pp. 162--172, 2011.

\bibitem{2011_CDC_CSPB}
M.~Chertkov, M.~Stepanov, F.~Pan, and R.~Baldick, ``Exact and efficient
  algorithm to discover extreme stochastic events in wind generation over
  transmission power grids,'' in \emph{2011 50th {IEEE} Conference on Decision
  and Control and European Control Conference ({CDC}-{ECC})}, 2011, pp.
  2174--2180.

\bibitem{1993_IEEETPS_TaskForce}
W.~W. Price, H.-D. Chiang, H.~K. Clark, C.~Concordia, D.~C. Lee, J.~C. Hsu,
  S.~Ihara, C.~A. King, C.~J. Lin, Y.~Mansour, K.~Srinivasan, C.~W. Taylor, and
  E.~Vaahedi, ``Load representation for dynamic performance analysys,''
  \emph{IEEE Trans. Power Syst.}, vol.~8, no.~2, pp. 472--482, 1993.

\bibitem{1973_PRA_MSR}
P.~C. Martin, E.~D. Siggia, and H.~A. Rose, ``Statistical dynamics of classical
  systems,'' \emph{Phys. Rev. A}, vol.~8, no.~1, pp. 423--437, 1973.

\bibitem{Pontryagin}
L.~S. Pontryagin, V.~G. Boltyanskii, R.~V. Gamkrelidze, and E.~F. Mishchenko,
  \emph{The mathematical theory of optimal processes}.\hskip 1em plus 0.5em
  minus 0.4em\relax New York: Interscience Publishers, 1962.

\bibitem{Betts}
J.~T. Betts, \emph{Practical methods for optimal control using nonlinear
  pro\-gram\-ming}.\hskip 1em plus 0.5em minus 0.4em\relax Philadelphia: SIAM,
  2001.

\bibitem{Bellman}
R.~Bellman, \emph{Dynamic programming}.\hskip 1em plus 0.5em minus 0.4em\relax
  Princeton: Princeton University Press, 1957.

\bibitem{2004_PD_FKLMT}
G.~Falkovich, I.~Kolokolov, V.~Lebedev, V.~Mezentsev, and S.~Turitsyn,
  ``Non-{G}aussian error probability in optical soliton transmission,''
  \emph{Physica D}, vol. 195, no. 1--2, pp. 1--28, 2004.

\end{thebibliography}

\end{document}